\documentclass{article}

\usepackage{klm_makros}

\begin{document}

\title{Box complexes, neighborhood complexes,\\ and the chromatic number}

\author{
{\sc P\'eter Csorba\thanks{Supported by the joint Berlin/Z\"urich graduate
                           program ``Combinatorics, Geometry, and Computation
                           (CGC),'' financed by ETH Z\"urich and the Deutsche
                           Forschungsgemeinschaft (DFG grant GRK 588/1).}
}\\
{\footnotesize  Institute of Theoretical Computer Science}\\[-1.5mm]
{\footnotesize  ETH Z\"urich, 8092 Z\"urich, Switzerland}\\[-1.5mm]
{\footnotesize E-mail: \texttt{pcsorba@inf.ethz.ch}}
\and \addtocounter{footnote}{5}
{\sc Carsten Lange\thanks{Supported by the DFG Sonderforschungsbereich 288
                          ``Differentialgeometrie und Quantenphysik'' in Berlin.}
}\\
{\footnotesize Institute of Mathematics, MA 6-2}\\[-1.5mm]
{\footnotesize TU Berlin D-10623 Berlin, Germany}\\[-1.5mm]
{\footnotesize  E-mail: \texttt{lange@math.tu-berlin.de}}
\and
{\sc Ingo Schurr$^{*}$}\\
{\footnotesize  Institute of Theoretical Computer Science}\\[-1.5mm]
{\footnotesize  ETH Z\"urich, 8092 Z\"urich, Switzerland}\\[-1.5mm]
{\footnotesize E-mail: \texttt{schurr@inf.ethz.ch}}
\and
{\sc Arnold Wa{\ss}mer$^{*}$}\\
{\footnotesize Institute of Mathematics, MA 6-2}\\[-1.5mm]
{\footnotesize TU Berlin D-10623 Berlin, Germany}\\[-1.5mm]
{\footnotesize  E-mail: \texttt{wassmer@math.tu-berlin.de}}
}

\date{October 21, 2003}

\maketitle

\begin{abstract}
Lov\'asz's striking proof of Kneser's conjecture from 1978 using the Borsuk--Ulam theorem
provides a lower bound on the chromatic number $\chr G$ of a graph $G$.
We introduce the \emph{shore subdivision} of simplicial complexes and use it to show
an upper bound to this topological lower bound and to construct a strong $\Z_{2}$-deformation retraction
from the box complex (in the version introduced by Matou{\v s}ek and Ziegler) to the Lov\'asz complex.
In the process, we analyze and clarify the combinatorics of the complexes involved and
link their structure via several ``intermediate'' complexes.
\end{abstract}
\section{Introduction}

The topological method in graph theory was introduced by Lov\'asz \cite{L} to prove Kneser's conjecture \cite{Kneser}. 
The pattern to obtain a lower bound of the chromatic number $\chr G$ of a graph $G$ is to associate a topological space 
and bound the chromatic number by a topological invariant of this space, e.g. connectivity or $\Z_{2}$-index. In 
this note we present a subdivision technique that shows that the complex $\lovaszcomplex G$ which Lov{\'a}sz used (and 
which we call Lov{\'a}sz complex for that reason) is a $\Z_{2}$-deformation retract of the box complex $\boxcomplex G$ 
described by Matou{\v s}ek and Ziegler~\cite{MZ}. The advantage of the box complex is that for any graph homomorphism 
$f: G \longrightarrow H$ one obtains an induced simplicial $\Z_{2}$-map $\boxcomplex f : \boxcomplex G \longrightarrow \boxcomplex H$. This functorial 
property gives elegant conceptual proofs which was not the case for the Lov{\'a}sz complex. Walker~\cite{W} constructed a 
$\Z_{2}$-map $\varphi: \| \lovaszcomplex G \| \longrightarrow \|\lovaszcomplex H \|$. Such a map could also be constructed 
using $\boxcomplex f$ and the $\Z_{2}$-deformation retraction constructed below.

The box complex of a graph yields a lower bound for its chromatic number: $\indx {\boxcomplex G} + 2 \leq \chr G$. 
It is known that this topological bound can get arbitrarily bad: Walker~\cite{W} shows that if a graph $G$ does not 
contain a $\KB 22$ then the associated invariant yields $3$ as largest possible lower bound for the chromatic number 
$\chr G$. In section~\ref{arnold} we generalize this result to the following statement: If $G$ does not contain a complete bipartite 
graph $\KB {\ell}m$ then the index of the box complex $\boxcomplex G$ is bounded by $\ell +m-3$ and this bound is sharp. 

Finally, we show in section~\ref{sec:peter} that $\lovaszcomplex G$ is $\Z_{2}$-isomorphic to a subcomplex of the shore 
subdivision of the box complex $\boxcomplex G$ (which is introduced in section~\ref{shoresubdivision}) and that this 
copy of $\lovaszcomplex G$ is a strong $\Z_{2}$-deformation retract of ${\boxcomplex G}$.

\section{Preliminaries}
In this section we recall some basic facts of graphs and simplicial complexes to fix notation. The interested reader 
is referred to~\cite{M} or~\cite{B} for details.

\subpart{Graphs} 
Any graph $G$ considered will be assumed to be finite, simple, connected, and undirected, i.e. $G$ is given by a finite 
set $\vertices G$ of \emph{nodes} (we use \emph{vertices} for associated complexes) and a set of \emph{edges} 
$\edges G \subseteq \binom{\vertices G}{2}$. A \emph{proper graph coloring with $n$ colors} is a homomorphism 
$c: G\to\K n$, where $\K n$ is the complete graph on $n$ nodes and the \emph{chromatic number} $\chr G$ of $G$ is the 
smallest $n$ such that there exists a proper graph coloring of $G$ with $n$ colors. The {\em neighborhood} $\neighbors u$ 
of $u \in \vertices G$ is the set of all nodes adjacent to $u$. For a set of nodes $A \subseteq \vertices G$ a node $v$ 
is in the {\em common neighborhood} $\cneighbors A$ of $A$, if $v$ is adjacent to all $a\in A$; we define 
$\cneighbors \emptyset := \vertices G$. For $A \subseteq B \subseteq \vertices G$ the common neighborhood relation satisfies 
(a) $A \cap \cneighbors A = \emptyset$, 
(b) $\cneighbors B \subseteq \cneighbors A$, 
(c) $A \subseteq \cncn A$, and
(d) $\cneighbors A = \cncncn {A}$. 
For two disjoint sets of nodes $A,B \subseteq \vertices G$ we define $\bipsubgraph G{A}{B}$ as the (not necessarily induced) 
subgraph of $G$ with node set $\vertices { \bipsubgraph G{A}{B} } = A \cup B$ and all edges $\{ a,b \} \in \edges G$ with 
$a\in A$ and $b\in B$. In this notation $\cneighbors A$ is the inclusion-maximal set $B$ such that $\bipsubgraph  G AB$ is 
complete bipartite.

\subpart{Simplicial Complexes}
An \emph{abstract simplicial complex} $\complex K$ is a finite hereditary set system.
We denote its vertex set by $\vertices {\complex K}$ and its barycentric subdivision by $\sd {\complex K}$.
For sets $A,B$ define $\pair AB := \set{(a,0)}{a\in A} \cup \set{(b,1)}{b\in B}$.
An important construction in the category of simplicial complexes is the {\em join operation}. 
For two simplicial complexes $\complex K$ and $\complex L$ the join ${\complex K}*{\complex L}$ is defined as
$\{ \pair FG \, | \, F\in \complex K \text{ and } G\in \complex L \}$. Any abstract simplicial
complex $\complex K$ can be realized as a topological space $\| \complex K \|$ in $\R^{d}$ for some $d$.

\subpart{$\Z_{2}$-spaces}
A \emph{$\Z_2$-space} is a topological space $X$ together with a homeomorphism $\nu: X \to X$ that is self-inverse 
and free, i.e. has no fixed points. The map $\nu$ is called \emph{free $\Z_{2}$-action}. The fundamental example 
for a $\Z_{2}$-space is the $d$-sphere $S^{d}$ together with the antipodal map $\nu(x)=-x$. A continuous map $f$ 
between $\Z_{2}$-spaces $(X,\nu)$ and $(Y,\mu)$ is \emph{$\Z_{2}$-equivariant} (or a \emph{$\Z_{2}$-map} for 
simplicity) if $f$ commutes with the $\Z_{2}$-actions, i.e. $f \circ \nu = \mu \circ f$. A simplicial complex 
$(\complex K,\nu)$ is a \emph{simplicial $\Z_{2}$-space} if $\nu: \complex K \to \complex K$ is a simplicial map 
such that $\geom \nu $ is a free $\Z_{2}$-action on $\geom{ \complex K }$. A \emph{simplicial $\Z_{2}$-equivariant} 
map $f$ is a simplicial map between two simplicial $\Z_{2}$-spaces that commutes with the simplicial $\Z_{2}$-actions.  
The {\em index} of a $\Z_2$-space $(X,\nu)$ is the smallest $d$ such that there is a $\Z_{2}$-map $f:X \to S^{d}$,
i.e. $f\circ \nu = - f$. The Borsuk--Ulam theorem provides the index for spheres: $\indx{ S^{d}} = d$. Since the 
$\Z_2$-actions are usually clear, we tend to refer to a $\Z_2$-space $\complex K$ without explicit reference to $\nu$.

\subpart{Chain Notation}
We denote by $\A$ a chain $A_{1} \subset \ldots \subset A_{p}$ of subsets of $\vertices G$. A chain $\A$ will be of 
length $p$ and a chain $\B$ of length $q$. For $1 \leq t \leq p$ we denote by $\A_{\leq t}$ the chain 
$A_{1} \subset \ldots \subset A_{t}$. A similar convention is used for $\A_{\geq t}$. For chains $\A$, $\B$ satisfying 
$A_{p} \subseteq B_{1}$ the chain $A_{1} \subset \ldots \subset A_{p}\subseteq B_{1} \subset \ldots \subset B_{q}$ 
will be denoted by $\A \sqsubseteq \B$, where we omit $A_{p}$ or $B_{1}$ in case $A_{p}=B_{1}$. 
If a map $f$ preserves (resp. reverses) orders, we write $f(\A)$ instead 
of $f(A_{1}) \subseteq \ldots \subseteq f(A_{p})$ (resp. $f(A_{p}) \subseteq \ldots \subseteq f(A_{1})$).

\subpart{Neighborhood Complex}
The \emph{neighborhood complex} $\neighborhoodcomplex G$ of a graph $G$ has $\vertices G$ as vertices and the sets 
$A \subseteq \vertices G$ with $\cneighbors A \neq \emptyset$ as simplices.

\subpart{Lov\'asz Complex}
In general $\neighborhoodcomplex G$ is not a $\Z_2$-space. However, the neighborhood complex can be retracted to a 
$\Z_{2}$-subspace, the \emph{Lov\'asz complex}. This complex $\lovaszcomplex G$ is the subcomplex of
$\sd {\neighborhoodcomplex G}$ induced by the vertices that are fixed points of $\text{CN}^{2}$. The Lov\'asz complex is
\[
  \lovaszcomplex G = \set{ \A }
                         { \A \text{ a chain of node sets of $G$ with } \A=\cncn{\A}  }
\]
which is a $\Z_{2}$-space with $\Z_2$-action $\text{CN}$.

\subpart{Box Complex}
Different versions of a box complex are described by Alon, Frankl, and Lov{\'a}sz~\cite{AFL},  Sarkaria~\cite{S}, 
K{\v r}{\' i}{\v z}~\cite{Kriz}, and Matou{\v s}ek and Ziegler~\cite{MZ}. The \emph{box complex} $\boxcomplex G$ 
of $G$ in which we are interested is the one introduced by Matou{\v s}ek and Ziegler and is defined by
\begin{align*}
\boxcomplex G  :&= \set{ \pair AB }
                       { A,B \in \neighborhoodcomplex G \text{ and }
                         \bipsubgraph{G}{A}{B} \text{ is complete bipartite}
                       }\\                        
                &= \set{ \pair AB }
                       { A,B \in \neighborhoodcomplex G,\text{ }
                         A \subseteq \cneighbors B, \text{ and } B \subseteq \cneighbors A }.     
\end{align*}
The vertices of the box complex are $V_{1} := \pair {\{ v \}}\emptyset$ and $V_{2} := \pair \emptyset{\{ v \}}$ for 
all vertices of $G$. The subcomplexes of $\boxcomplex G$ induced by $V_{1}$ and $V_{2}$ are disjoint  subcomplexes 
of $\boxcomplex G$ that are both isomorphic to the neighborhood complex $\neighborhoodcomplex G$. We refer to these 
two copies as \emph{shores} of the box complex. The box complex is endowed with a $\Z_{2}$-action $\nu$ which 
interchanges the shores.

\section{Shore Subdivision and Useful Subcomplexes}
\label{shoresubdivision}
\subpart{Shore Subdivision}
More general, for a simplicial complex $\complex K$ and any partition $V_1 \sqcup V_2$ of its vertex set, we call 
the simplicial subcomplexes $\complex K_{1}$ and $\complex K_{2}$ induced by the vertex sets $V_{1}$ and $V_{2}$ 
its \emph{shores}. The \emph{shore subdivision} of $\complex K$ is
   \[
    \ssd{\complex K} := \set{\sd {\sigma \cap \complex K_{1}} * \sd {\sigma \cap \complex K_2}}
                            {\sigma \in \complex K}.
   \]
The shores of the box complex define a partition of the vertex set which 
allows us to define the shore subdivision $\ssd {\boxcomplex G}$ of the box complex $\boxcomplex G$. The vertices 
of $\ssd {\boxcomplex {G}}$ are of type $\pair A\emptyset$ and  $\pair \emptyset A$ where $\emptyset \neq A \subset V(G)$ 
with $\cneighbors A \neq \emptyset$. A simplex of $\ssd {\boxcomplex {G}}$ is denoted by $\pair \A\B$ (the simplex 
spanned by the vertices $\pair A\emptyset$ and $\pair \emptyset B$ where $A\in\A$, $B\in\B$).

\subpart{Doubled Lov{\' a}sz Complex}
The map $\scncn: \ssd{ \boxcomplex G} \to \ssd{ \boxcomplex G}$ defined on the vertices by
$\pair {A}{\emptyset}   \mapsto  \pair {\cncn A}{\emptyset}$ and 
$\pair {\emptyset}{A}   \mapsto  \pair {\emptyset}{\cncn A}$
is simplicial and $\Z_{2}$-equivariant. We refer to its image $\im \scncn$ as \emph{doubled Lov\'asz complex} 
$\dlovaszcomplex G$. It is
\[
  \dlovaszcomplex G = \set{ \pair \A\B  }
                          { \begin{matrix}
                                \A, \B \in \lovaszcomplex G,\\
				\bipsubgraph{G}{A}{B} \text{ is complete bipartite for all }A\in \A, B \in \B
                             \end{matrix}
			  }.
\]

\noindent
A copy of the Lov{\'a}sz complex can be found on each shore of $\dlovaszcomplex G \subset \ssd {\boxcomplex G}$, 
but these copies do not respect the induced $\Z_{2}$-action. 

\subpart{Halved Doubled Lov{\' a}sz Complex}
We partition the vertex set of the doubled Lov\'asz complex $\dlovaszcomplex G$ into pairs of type 
$\{ \pair A\emptyset ,\pair \emptyset{\cneighbors A} \}$ to define a simplicial $\Z_{2}$-map 
$\jumpf : \dlovaszcomplex G \rightarrow \dlovaszcomplex G$.
Our aim is to specify one vertex for every pair and map both vertices of a pair to this chosen ``smaller'' vertex. 
To do this we refine the partial order by cardinality to a linear order ``$\prec$'' on the vertices of the original 
Lov\'asz complex $\lovaszcomplex G$ using the lexicographic order:
\[
  A  \prec B 
  \quad :\Longleftrightarrow \quad 
  \begin{cases}
     |A| < |B| \text{ or } \\
     |A| = |B| \text{ and } A <_{\text{lex}} B.
  \end{cases}
\]
In fact any refinement would work in the following. A partial order on the vertices of the doubled 
Lov\'asz complex $ \dlovaszcomplex G $ is now obtained:
\[
  \pair A\emptyset \prec \pair {\emptyset}{\cneighbors A} 
  \quad :\Longleftrightarrow \quad
  A \prec \cneighbors A.
\]
We define the map $\jumpf$ using this partial order by % on the vertices: 
$\jump{\pair {A}{\emptyset}}  :=  \min_\prec \{\pair {A}{\emptyset} , \pair {\emptyset}{\cneighbors A} \}$
and $\jump{\pair {\emptyset}{B}} :=  \min_\prec \{\pair {\emptyset}{B} , \pair {\cneighbors B}{\emptyset} \}$.
Since the image $\im \jumpf$ has half as many vertices as $\dlovaszcomplex G$, we refer to 
$\im \jumpf$ as \emph{halved doubled Lov\'asz complex} $\hdlovaszcomplex G$. 

\begin{figure}[t]
   (a)\includegraphics[width=3.5cm]{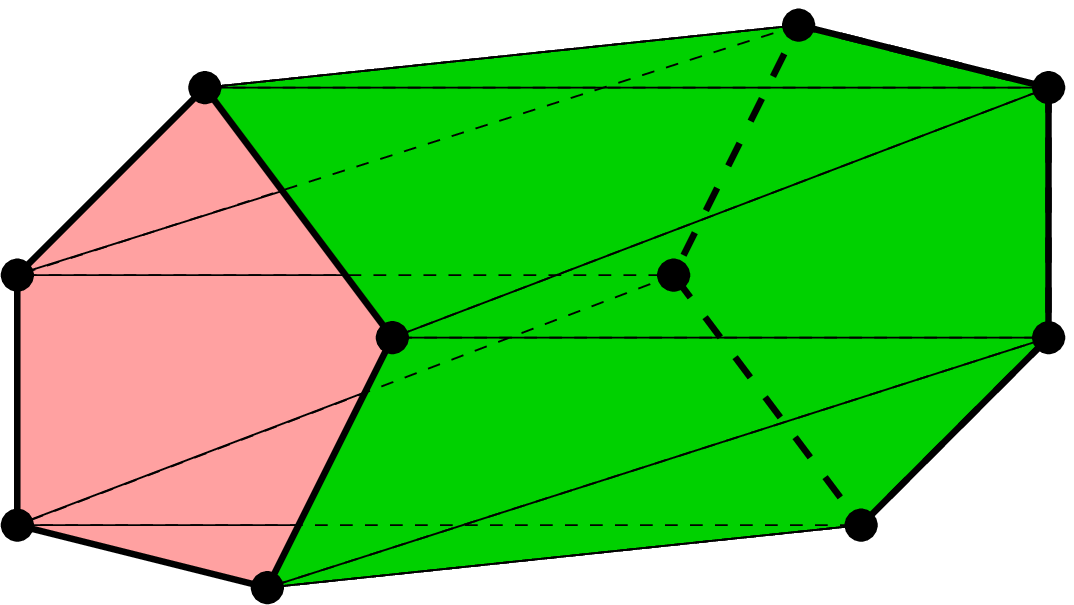}
   (b)\includegraphics[width=3.5cm]{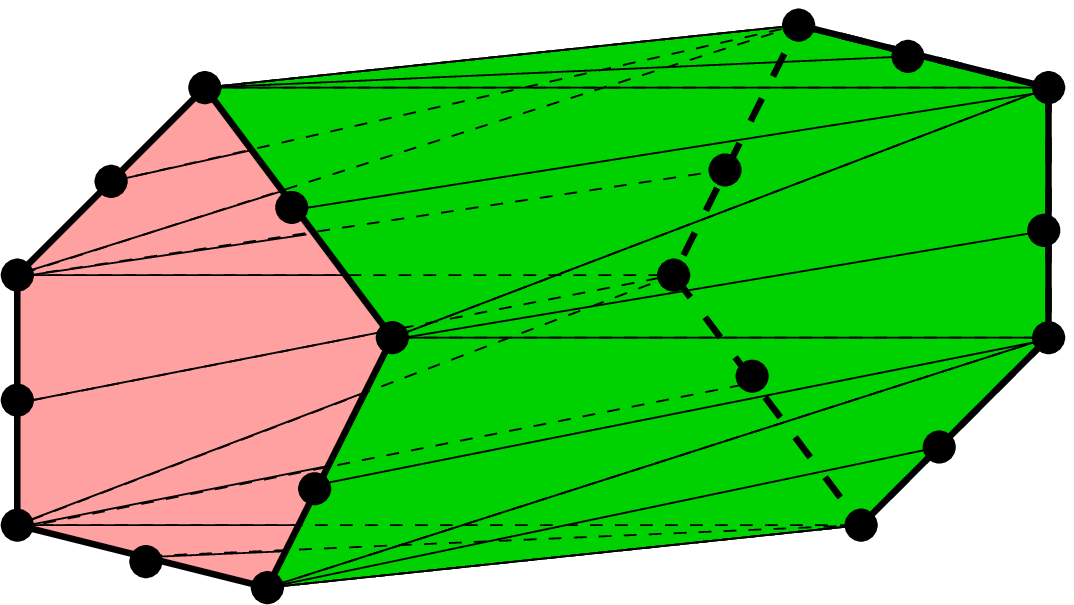}
   (c)\includegraphics[width=3.5cm]{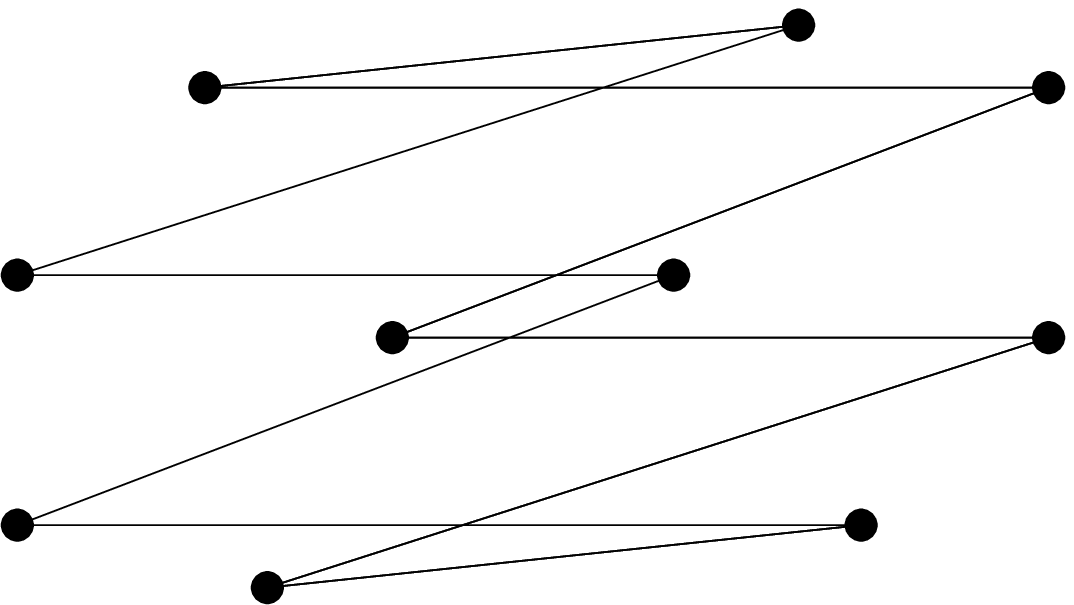}
   \caption{(a) $\boxcomplex {C_{5}}$; 
            (b) $\ssd {\boxcomplex {C_{5}}} = \dlovaszcomplex {C_{5}}$; and 
	    (c) $\hdlovaszcomplex {C_{5}}$}  
   \label{illustration}
\end{figure}

\subpart{An example}
The neighborhood complex $\neighborhoodcomplex {C_{5}}$ of the $5$-cycle $C_{5}$ is the $5$-cycle; its Lov{\' a}sz 
complex $\lovaszcomplex {C_{5}}$ is the $10$-cycle $C_{10}$. The box complex $\boxcomplex {C_{5}}$ consists of two 
copies of $\neighborhoodcomplex {C_{5}}$ (the two shores) such that simplices of different shores are joined iff 
their vertex sets are common neighbors of each other. The shore subdivision $\ssd {\boxcomplex {C_{5}}}$ is a 
subdivision of the box complex induced from a barycentric subdivision of the shores. The map $\scncn$ maps a vertex 
of $\ssd {\boxcomplex {C_{5}}}$ to the common neighborhood of its common neighborhood. In our example, every vertex 
is mapped to itself, hence $\ssd {\boxcomplex {C_{5}}} = \dlovaszcomplex {C_{5}}$. The partitioning of the vertex set 
of $\dlovaszcomplex {C_{5}}$ into pairs of type $(\pair A\emptyset , \pair \emptyset{\cneighbors A})$ can be visualized 
by edges of $\dlovaszcomplex {C_{5}}$ that connect singletons from one shore with two-element sets from the other. 
The smaller vertex of each such pair is actually a vertex of the original box complex $\boxcomplex {C_{5}}$.
Hence the map $\jumpf$ collapses all edges of type $(\pair A\emptyset, \pair \emptyset{\cneighbors A})$, which yields
the halved doubled Lov{\'a}sz complex $\hdlovaszcomplex G$. The maps $f_{i}$ introduced in section~\ref{sec:peter} are
these collapses and they are used to show that $\lovaszcomplex G$ is a $\Z_{2}$-deformation retract of 
$\ssd{\boxcomplex G}$. All these complexes are illustrated in Figure~\ref{illustration}. 

\section{The $\mathbf{K_{l,m} }$-Theorem   }
\label{arnold}

\begin{theorem} \label{thm:upper_bound}
  If a graph $G$ does not contain a complete bipartite subgraph
  $\KB {\ell}{m}$ then the index of its box complex is bounded by
  \[
    \indx{\boxcomplex G } \leq \ell + m - 3.
  \]
\end{theorem}
Since $\indx {\boxcomplex {\K {\ell + m - 1}}} = \ell + m - 3$, the statement of the theorem is best possible. On the 
other hand, we obtain $\indx {\boxcomplex {\K {k,k}}} \leq k-1$, but it can be shown that 
$\indx {\boxcomplex {\K {k,k}}} = 0$. So the gap in the inequality can arbitrarily large.

We give two proofs for this theorem. The first one uses the shore subdivision and the halved doubled Lov{\'a}sz complex,
the other is a direct argument on $\lovaszcomplex G$ along the lines of Walker~\cite{W}.

\smallskip
\begin{proof}(using Shore Subdivision)
Let $\Phi: \ssd{\boxcomplex G} \to \ssd{\boxcomplex G}$ be the simplicial $\Z_{2}$-map defined by $\jumpf \circ \scncn$. 
Using that the index is dominated by dimension, it suffices to show the last inequality of
\[
  \indx {\boxcomplex G} =\indx {\ssd {\boxcomplex G}} \leq \indx {\im \Phi} \leq \Dim {\im \Phi} \leq  \ell + m - 3.
\]
\noindent
To estimate the dimension of $\im \Phi = \hdlovaszcomplex G$, we use that the graph $G$ does not contain a 
$\KB{\ell}{m} $ as a subgraph and  assume without loss of generality that $\ell \leq m$.  
A vertex of $\hdlovaszcomplex G$ or $\dlovaszcomplex G$ of the form $\pair A\emptyset$ or $\pair \emptyset{A}$
is called \emph{small} if $|A| <  \ell$, \emph{medium} if $\ell  \leq |A| <  m$, and \emph{large} if $m  \leq |A|$.
For $\ell=m$ there are no medium vertices.
Let $\sigma = \pair \A\B$ be a simplex of $\hdlovaszcomplex G$ and consider  the set of vertices
   \[
     M_{\sigma} := \jumpf^{-1}(\sigma) =
     \bigcup_{A \in \A} \{ \pair {A}\emptyset,\pair \emptyset{\cneighbors {A}} \}
     \cup
     \bigcup_{B \in \B} \{ \pair {\cneighbors {B}}\emptyset,\pair \emptyset{B} \}.
   \]
Clearly, $|M_{\sigma}|$ is at most twice $| V(\sigma) |$.
If $\sigma$ has a large vertex $\pair A\emptyset$, then the vertex 
$\pair \emptyset{\cneighbors A}$ must be small, otherwise $G$ would contain a $\KB{\ell}{m}$. 
Hence there are at most $2\cdot 2(\ell - 1)$  many vertices in $M_{\sigma}$ that are large or small.
Since the number of medium vertices is at most $2(m-\ell)$, we have 
   \[ 
     | M_{\sigma}| \leq 2\cdot 2(\ell - 1)+2(m-\ell )=2(\ell +m-2).
   \]
Hence $| \vertices \sigma | \leq \ell +m-2$ for all $\sigma$, and therefore $\Dim {\hdlovaszcomplex G}$ is at most $\ell +m-3$.
\end{proof}

\smallskip

\smallskip
\begin{proof}(using Lov{\'a}sz Complex)
It suffices to prove $\Dim {\lovaszcomplex G } \le \ell + m - 3$ since 
$\indx{\boxcomplex G } = \indx{\lovaszcomplex G } \le \Dim {\lovaszcomplex G }$, (\cite{MZ} or section~\ref{sec:peter}).
Without loss of generality let $\ell \leq m$ and consider a simplex $\sigma = A_{1} \subset \ldots \subset A_{p}$ of 
$\lovaszcomplex G$ of maximal dimension $p-1$. If $p < \ell$ we are done. Suppose therefore that $p \geq \ell$. Then 
$\bipsubgraph G{A_{\ell}}{\cneighbors {A_{\ell}} }$ is a bipartite subgraph of $G$ and we have $| A_{\ell} | \geq \ell$ 
as well as $| \cneighbors {A_{\ell}} | \geq p - \ell +1$. The assumption that $G$ does not contain a $K_{\ell ,m}$ 
implies that $m > p - \ell + 1$, i.e.  $\Dim \sigma \leq \ell +m - 3$. 
\end{proof}

%%%%%%%%%%%%%%%%%%%%%%%%%%%%%%%%%%%%%%%%%%%%%%%%%%%%%%%%%%%%%%%%%%
%%
%%                      This is the LaTeX2e file for
%%
%%                    L(G) is Z_2-def., retract of B(G)
%%
%%                                  by
%%
%%%%%%%%%%%%%%%%%%%%%%%%%%%%%%%%%%%%%%%%%%%%%%%%%%%%%%%%%%%%%%%%%%%

\section{$\lovaszcomplex G$ as a $\Z_{2}$-Deformation Retract of $\boxcomplex G$}
\label{sec:peter}

\begin{theorem}
   The Lov\'asz complex $\lovaszcomplex G$ and the halved doubled Lov{\'a}sz complex $\hdlovaszcomplex G$ are $\Z_2$-isomorphic.
\end{theorem}

\medskip
\begin{proof}
   First we have $|\vertices {\lovaszcomplex G}| = |\vertices {\hdlovaszcomplex G}|$ since each shore of $\dlovaszcomplex G$ 
   is isomorphic (but not $\Z_{2}$-isomorphic) to $\lovaszcomplex G$. To define a simplicial $\Z_{2}$-map 
   $f:\lovaszcomplex G \rightarrow \hdlovaszcomplex G$, we partition $\vertices {\lovaszcomplex G}$ into
   \begin{align*}
      S :=\set{A }{\begin{matrix}
                      A \in \vertices {\lovaszcomplex G} \text{ and } \\
                      \jump{A \uplus\emptyset}=A\uplus\emptyset
                   \end{matrix}}
      \quad \text{and}\quad
      J :=\set{A }{\begin{matrix}
                      A \in \vertices{\lovaszcomplex G} \text{ and } \\
                      \jump{A \uplus\emptyset}=\emptyset\uplus\cneighbors A
                   \end{matrix}},
   \end{align*}
   (where ``$S$'' and ``$J$'' denote the vertices that stay fixed or jump to their neighbor), and set
   \[
     f(A) := \begin{cases}
                A \uplus\emptyset            & \text{ if } A\in S\\
                \emptyset\uplus\cneighbors A & \text{ if } A\in J.
             \end{cases}
   \]
   This map is a bijection between the vertex sets, surjective, simplicial, and $\Z_2$-~equivariant. For simpliciality, 
   consider a simplex  $\A$ in $\lovaszcomplex G$. Let $t$ denote the largest index $i$ such that $A_{i}$ is mapped onto 
   the first shore. The image of $\A$ under $f$ is $\pair {\A_{\leq t}}{\cneighbors{\A_{\geq t+1}}}$. This is a simplex 
   since $\bipsubgraph{G}{A_t}{\cneighbors{A_{t+1}}}$ is complete bipartite. For surjectivity consider a simplex  
   $\pair {\A}{\B}$ of $\hdlovaszcomplex G$, i.e. $\bipsubgraph{G}{A_p}{B_q}$ is complete bipartite. This is the image 
   of the simplex $\A \subseteq \cneighbors{\B}$ of $\lovaszcomplex G$.
\end{proof}

\begin{theorem}
   The halved doubled Lov{\'a}sz complex $\hdlovaszcomplex G$ is a strong $\Z_2$-defor\-mation retract of the box complex 
   $\boxcomplex G$.
\end{theorem}

\begin{proof}
   First we observe that $\| \dlovaszcomplex G \|$ is a strong $\Z_2$-deformation retract of $\| \ssd{\boxcomplex G} \|$. 
   This follows from the fact that a closure operator induces a strong deformation retraction from its domain to its image 
   (\cite{B}, \cite{M}). Explicitly, this map is obtained by sending each point $p \in \| \ssd{\boxcomplex G} \|$ towards 
   $\| \cncnf \| (p)$ with uniform speed, which is $\Z_{2}$-equivariant at any time of the deformation.

   To show that $\| \hdlovaszcomplex G \|$ is a strong $\Z_2$-deformation retract of $\| \dlovaszcomplex G \|$, we define 
   simplicial complexes and simplicial $\Z_{2}$-maps
   \[
     \dlovaszcomplex G =: S_{0}
              \stackrel{f_0}{\longrightarrow} S_1
              \stackrel{f_1}{\longrightarrow} \ldots
              \stackrel{f_{N}}{\longrightarrow} S_{N+1} := \hdlovaszcomplex G
   \]
   such that $S_{i+1}$ is a $\Z_2$-subcomplex of $S_{i}$ and $S_{i+1}$ is a strong $\Z_2$-deformation retract of 
   $S_{i}$. The composition of the $f_i$ yields the earlier defined map $\jumpf$, i.e. 
   $\jumpf=f_{N}\circ \dots \circ f_1 \circ f_0$. To construct $S_{i+1}$ inductively from $S_{i}$, 
   we consider $X := \max_{\prec} \set{Y \in J}{\pair Y\emptyset \in S_{i}}$ and obtain $S_{i+1}$  from $S_{i}$ by deleting 
   each simplex of $S_{i}$ that contains $X\uplus\emptyset$ or its $\Z_2$-pair $\emptyset\uplus X$, i.e.
   \[
     S_{i+1} := \set{\sigma}
                    {\sigma\in S_i \text{ and } X\uplus\emptyset \not\in \sigma \text{ and } \emptyset\uplus X \not\in \sigma }.
   \]
   The maximality of $X$ implies that a maximal simplex which contains $\pair X\emptyset$ (resp. $\pair \emptyset X$) does also 
   contain $\pair \emptyset{\cneighbors X}$ (resp. $\pair {\cneighbors X}\emptyset$). Hence the map $f_i$ defined on the vertices 
   $v \in \vertices {S_{i}}$ via
   \[
     f_{i}(v) := \begin{cases}
                    \emptyset\uplus\cneighbors X & \text{if } v=X \uplus\emptyset \\
                    \cneighbors X\uplus\emptyset & \text{if } v=\emptyset\uplus X \\
                    v                            & \text{otherwise}
                 \end{cases}
   \]
   \noindent
   is simplicial and $\Z_{2}$-equivariant.

   Thus $F\colon \| S_{i} \| \times [0,1] \to \| S_{i} \|$ given by $F(x,t):=t\cdot x+(1-t) \cdot \|f_i \| (x)$ 
   is a well-defined $\Z_2$-homotopy from $\| f_i \|$ to $\id_{\| S_{i} \|}$ that fixes $\| S_{i+1} \|$.
\end{proof}

\subsection*{Acknowledgments}

The authors thank G\"unter M.\ Ziegler for bringing the problems
studied in this paper to their attention, and Tibor Szab\'o for numerous discussions.

\end{document}